\documentclass[reqno,onecolumn,12pt]{amsart}
\usepackage{amsfonts}
\usepackage{amssymb}
\usepackage{amsmath}
\usepackage{graphicx}
\usepackage{amscd}

\newtheorem{theorem}{Theorem}
\theoremstyle{plain}

\newtheorem{corollary}{Corollary}

\newtheorem{definition}{Definition}

\numberwithin{equation}{section}

\input{tcilatex}

\begin{document}
\author{}
\title{}
\maketitle

\begin{center}
\thispagestyle{empty}\textbf{ON THE IDENTITIES INVOLVING SPECIAL POLYNOMIALS
ARISING FROM POINT OF VIEW OF FRACTIONAL CALCULUS}

\bigskip \textbf{By}

\hspace{0.05cm}

\hspace{0.05cm}\textbf{Serkan Araci}$^{1}$\textbf{, Erdo\u{g}an \c{S}en}$%
^{2,3}$\textbf{, Mehmet Acikgoz}$^{1}$\textbf{\ and Kamil Oru\c{c}o\u{g}lu}$%
^{3}$

\hspace{-0.5cm}

\bigskip $^{1}$University of Gaziantep, Faculty of Arts and Science,
Department of Mathematics, 27310 Gaziantep TURKEY

\hspace{0.5cm}

$^{2}$Department of Mathematics, Faculty of Science and Letters, Namik Kemal
University, 59030 Tekirda\u{g}, TURKEY\bigskip

$^{3}$Department of Mathematics Engineering, Istanbul Technical University,
Maslak, 34469 Istanbul, TURKEY\bigskip \hspace{0.05cm}

mtsrkn@hotmail.com; erdogan.math@gmail.com; acikgoz@gantep.edu.tr;
koruc@itu.edu.tr

\hspace{-0.5cm}

\textbf{{\Large {Abstract}}}

\hspace{-0.5cm}
\end{center}

In the present paper, we investigate some interesting properties including
several special polynomials arising from Caputo-fractional derivative. From
our investigation, we derive a lot of interesting identities of several
special polynomials.

\bigskip


\hspace{-0.5cm}


\section{\textbf{Introduction}}

The concept of fractional is popularly appeared from a question raised in
the year 1695 by L'H\^{o}pital (1661-1704) to Leibniz (1646-1716), which
searched the way of Leibniz's notation $\frac{d^{n}y}{dx^{n}}$ for the
derivative of order $n\in 
\mathbb{N}
^{\ast }:=%
\mathbb{N}
\cup \left\{ 0\right\} $ when $n=\frac{1}{2}$ (What if $n=\frac{1}{2}$). In
his response, dated 30 September 1695, Leibniz wrote to L' H\^{o}pital as
follows:"$\cdots $This is an apparent paradox from which, one day, useful
consequences will be drawn$\cdots $" (see \cite{Podlubny}, \cite{Kilbas}).

Subsequent mention of fractional derivatives was made, in some context, by
Euler in 1730, Lagrange in 1772, Laplace in 1812, Lacroix in 1819, Fourier
in 1822, Liouville in 1832, Riemann in 1847, Greer in 1859, Holmgren in
1865, Gr\"{u}nwald in 1867, Letnikov in 1868, Sonin in 1869, Laurent in
1884, Nekrassov in 1888, Krug in 1890, and Weyl in 1917 (for details, see 
\cite{Podlubny}, \cite{Kilbas}).

One of the most recent works on the subject of fractional calculus is the
book of Podlubny \cite{Podlubny} published in 1999, which deals with theory
of fractional differential equations.

One of the fundamental functions of the fractional calculus is Euler's gamma
function $\Gamma \left( \xi \right) $, which generalizes the fractional "$n!$%
" defined by%
\begin{equation}
\Gamma \left( \xi \right) =\int_{0}^{\infty }t^{\xi -1}e^{-t}dt,
\label{equ.1}
\end{equation}%
which converges in the right half of the complex plane $\func{Re}\left( \xi
\right) >0$ (see \cite{Podlubny}, \cite{Kilbas}). This function satisfies
the following functional equations:%
\begin{equation}
\Gamma \left( \xi +1\right) =\xi \Gamma \left( \xi \right) \text{ \textit{and%
} }B\left( \xi ,\gamma \right) =\frac{\Gamma \left( \xi \right) \Gamma
\left( \gamma \right) }{\Gamma \left( \xi +\gamma \right) }  \label{equ.2}
\end{equation}%
where $B\left( \xi ,\gamma \right) $ is known as Beta function (see \cite%
{Podlubny}, \cite{Kilbas}, \cite{Kim 3}, \cite{Araci 1}).

The Mittag-Leffler function, which generalizes the exponential function $%
e^{z}$, is given by%
\begin{equation}
E_{\alpha }\left( z\right) =\sum_{n=0}^{\infty }\frac{z^{n}}{\Gamma \left(
\alpha n+1\right) }\text{ (see \cite{Podlubny}, \cite{Kilbas}).}
\label{equ.3}
\end{equation}

Obviously that%
\begin{equation*}
E_{1}\left( z\right) =e^{z}\text{.}
\end{equation*}

The Mittag-Leffler function plays a vital and important role in the concept
of \textit{fractional calculus}.

The generalization of Mittag-Leffler function is also defined by the
following series expansion:%
\begin{equation}
E_{\alpha ,\beta }\left( z\right) =\sum_{n=0}^{\infty }\frac{z^{n}}{\Gamma
\left( \alpha n+\beta \right) }\text{, }\left( \alpha >0,\beta >0\right) 
\text{ (see \cite{Podlubny}, \cite{Kilbas}).}  \label{equ.4}
\end{equation}

By (\ref{equ.4}), we have%
\begin{eqnarray}
E_{1,1}\left( z\right) &=&\sum_{k=0}^{\infty }\frac{z^{k}}{\Gamma \left(
k+1\right) }=e^{z}\text{,}  \notag \\
E_{1,2}\left( z\right) &=&\sum_{k=0}^{\infty }\frac{z^{k}}{\Gamma \left(
k+2\right) }=\frac{e^{z}-1}{z}\text{,}  \label{equ.5} \\
E_{1,3}\left( z\right) &=&\sum_{k=0}^{\infty }\frac{z^{k}}{\Gamma \left(
k+3\right) }=\frac{e^{z}-z-1}{z^{2}}\text{.}  \notag
\end{eqnarray}

Continuing this process, yields to%
\begin{equation*}
E_{1,m}\left( z\right) =\frac{1}{z^{m-1}}\left\{ e^{z}-\sum_{k=0}^{m-2}\frac{%
z^{k}}{k!}\right\} \text{ (see \cite{Podlubny}, \cite{Kilbas}).}
\end{equation*}

Note that variously generalizations of the Mittag-Leffler function are
studied by Humbert and Delerue \cite{Podlubny} and by Chak \cite{Chak}, were
further extended by Srivastava \cite{Kilbas}.

In \cite{Podlubny},\cite{Kilbas}, the Riemann-Liouville fractional integral
of order $\alpha $ for a function $f$ is defined by 
\begin{equation}
\mathcal{I}^{\left( \alpha \right) }f\left( t\right) =\frac{1}{\Gamma \left(
\alpha \right) }\int_{0}^{t}\left( t-s\right) ^{\alpha -1}f\left( s\right) ds%
\text{, }\left( f:\left( 0,\infty \right) \rightarrow 
\mathbb{R}
\text{ and }\alpha >0\right) \text{.}  \label{equ.6}
\end{equation}

The Caputo-fractional derivative of higher order for a continuos function $f$
is given by 
\begin{equation}
\mathcal{D}^{\left( \alpha \right) }f\left( t\right) =\frac{1}{\Gamma \left(
n-\alpha \right) }\int_{0}^{t}\frac{f^{\left( n\right) }\left( s\right) }{%
\left( t-s\right) ^{\alpha -n+1}}ds,\left( f:\left( 0,\infty \right)
\rightarrow 
\mathbb{R}
\text{ and }\alpha >0\right) ,\text{ }  \label{equ.7}
\end{equation}%
where $n$ is the smallest integer greater than or equal to $\alpha $ (see 
\cite{Podlubny}, \cite{Kilbas}).

From (\ref{equ.6}) and (\ref{equ.7}), we have%
\begin{equation}
\mathcal{D}^{\left( \alpha \right) }f\left( t\right) =\mathcal{D}^{\left(
k\right) }\left[ \mathcal{I}^{\left( k-\alpha \right) }f\left( t\right) %
\right] ,k\in 
\mathbb{N}
\text{.}  \label{equ.8}
\end{equation}

From (\ref{equ.8}), it follows that%
\begin{equation}
\mathcal{D}^{\left( \alpha \right) }t^{n}=\frac{\Gamma \left( n+1\right) }{%
\Gamma \left( n-\alpha +1\right) }t^{n-\alpha }\text{(see \cite{Podlubny}, 
\cite{Kilbas}).}  \label{equ.9}
\end{equation}

The fractional derivative of the product $fg$, which is called Leibniz rule,
is given by%
\begin{equation}
\mathcal{D}^{\left( \alpha \right) }\left[ f\left( t\right) g\left( t\right) %
\right] =\sum_{k=0}^{\infty }\binom{\alpha }{k}f^{\left( k\right) }\left(
t\right) \mathcal{D}^{\left( \alpha -k\right) }g\left( t\right) \text{ (see 
\cite{Podlubny}, \cite{Kilbas}).}  \label{equ.10}
\end{equation}

In the next section, we consider analogues of Bernoulli, Euler and Genocchi
polynomials which are derived from (\ref{equ.3}).

\section{\textbf{Analogues of Bernoulli, Euler and Genocchi polynomials and
their properties}}

Recently, analogues of Bernoulli, Euler and Genocchi polynomials were
studied by many mathematicians [4-24]. We are now ready to give the
definition of generating functions, corresponding to Mittag-Leffler
function, of Bernoulli, Euler and Genocchi polynomials.

\begin{definition}
\label{Definition 1}Let $\alpha >0$ and $\lambda >0$, define%
\begin{eqnarray*}
\mathcal{K}(x,z &:&\alpha \mid \lambda )=\sum_{n=0}^{\infty }\mathcal{B}%
_{n}\left( x:\alpha \mid \lambda \right) \frac{z^{n}}{n!}=\frac{z}{\lambda
E_{\alpha }\left( z\right) -1}e^{xz}\text{,} \\
\mathcal{I}(x,z &:&\alpha \mid \lambda )=\sum_{n=0}^{\infty }\mathcal{E}%
_{n}\left( x:\alpha \mid \lambda \right) \frac{z^{n}}{n!}=\frac{2}{\lambda
E_{\alpha }\left( z\right) +1}e^{xz}\text{,} \\
\mathcal{M}(x,z &:&\alpha \mid \lambda )=\sum_{n=0}^{\infty }\mathcal{G}%
_{n}\left( x:\alpha \mid \lambda \right) \frac{z^{n}}{n!}=\frac{2z}{\lambda
E_{\alpha }\left( z\right) +1}e^{xz}\text{,}
\end{eqnarray*}%
where $\mathcal{B}_{n}\left( x:\alpha \mid \lambda \right) $, $\mathcal{E}%
_{n}\left( x:\alpha \mid \lambda \right) $ and $\mathcal{G}_{n}\left(
x:\alpha \mid \lambda \right) $ are called, respectively, Bernoulli-type,
Euler-type and Genocchi-type polynomials.
\end{definition}

\begin{corollary}
Taking $\alpha =1$ in Definition \ref{Definition 1}, we have 
\begin{equation*}
\sum_{n=0}^{\infty }B_{n}\left( x\mid \lambda \right) \frac{z^{n}}{n!}=\frac{%
z}{\lambda e^{z}-1}e^{xz},\text{ }\sum_{n=0}^{\infty }E_{n}\left( x\mid
\lambda \right) \frac{z^{n}}{n!}=\frac{2}{\lambda e^{z}+1}e^{xz}\text{ and }%
\sum_{n=0}^{\infty }G_{n}\left( x\mid \lambda \right) \frac{z^{n}}{n!}=\frac{%
2z}{\lambda e^{z}+1}e^{xz}\text{,}
\end{equation*}%
where $B_{n}\left( x\mid \lambda \right) $, $E_{n}\left( x\mid \lambda
\right) $ and $G_{n}\left( x\mid \lambda \right) $ are called
Apostol-Bernoulli polynomials, Apostol-Euler polynomials and
Apostol-Genocchi polynomials, respectively (see \cite{Kim 5}, \cite{Sen}, 
\cite{Jolany}, \cite{Luo1}, \cite{He1}, \cite{He2}, \cite{Srivastava 1}).
\end{corollary}

\begin{corollary}
Substituting $\alpha =\lambda =1$ in Definition \ref{Definition 1}, we have 
\begin{equation*}
\sum_{n=0}^{\infty }B_{n}\left( x\right) \frac{z^{n}}{n!}=\frac{z}{e^{z}-1}%
e^{xz},\text{ }\sum_{n=0}^{\infty }E_{n}\left( x\right) \frac{z^{n}}{n!}=%
\frac{2}{e^{z}+1}e^{xz}\text{ and }\sum_{n=0}^{\infty }G_{n}\left( x\right) 
\frac{z^{n}}{n!}=\frac{2z}{e^{z}+1}e^{xz}\text{,}
\end{equation*}%
where $B_{n}\left( x\right) $, $E_{n}\left( x\right) $ and $G_{n}\left(
x\right) $ are called classical Bernoulli polynomials, classical Euler
polynomials and classical Genocchi polynomials, respectively (see \cite{Kim
5}, \cite{Sen}, \cite{Jolany}, \cite{Luo1}, \cite{He1}, \cite{He2}, \cite%
{Srivastava 1}).
\end{corollary}

Taking $x=0$ in the above definition, we have%
\begin{eqnarray*}
\mathcal{B}_{n}\left( 0:\alpha \mid \lambda \right) &:&=\underset{\text{%
Bernoulli-type number}}{\underbrace{\mathcal{B}_{n}\left( \alpha \mid
\lambda \right) }}\text{ } \\
\mathcal{E}_{n}\left( 0:\alpha \mid \lambda \right) &:&=\underset{\text{%
Euler-type number}}{\underbrace{\mathcal{E}_{n}\left( \alpha \mid \lambda
\right) }} \\
\text{ }\mathcal{G}_{n}\left( 0:\alpha \mid \lambda \right) &:&=\underset{%
\text{Genocchi-type number}}{\underbrace{\mathcal{G}_{n}\left( \alpha \mid
\lambda \right) }}
\end{eqnarray*}

and from the above, we write 
\begin{eqnarray}
\mathcal{K}(0,z &:&\alpha \mid \lambda ):=\mathcal{K}(z:\alpha \mid \lambda )%
\text{, }  \notag \\
\mathcal{I}(0,z &:&\alpha \mid \lambda ):=\mathcal{I}(z:\alpha \mid \lambda )%
\text{,\textit{\ }}  \label{equ.12} \\
\mathcal{M}(0,z &:&\alpha \mid \lambda ):=\mathcal{M}(z:\alpha \mid \lambda )%
\text{.}  \notag
\end{eqnarray}

Matching Definition \ref{Definition 1} and (\ref{equ.12}), we get the
following corollary.

\begin{corollary}
\label{Corollary 1}The following functional equations hold true:%
\begin{eqnarray*}
\mathcal{K}(x,z &:&\alpha \mid \lambda )=\mathcal{K}(z:\alpha \mid \lambda
)e^{xz}\text{,} \\
\mathcal{I}(x,z &:&\alpha \mid \lambda )=\mathcal{I}(z:\alpha \mid \lambda
)e^{xz}\text{,} \\
\mathcal{M}(x,z &:&\alpha \mid \lambda )=\mathcal{M}(z:\alpha \mid \lambda
)e^{xz}\text{.}
\end{eqnarray*}
\end{corollary}

By using (\ref{equ.4}) and Corollary \ref{Corollary 1}, becomes%
\begin{equation*}
\sum_{n=0}^{\infty }\mathcal{B}_{n}\left( x:\alpha \mid \lambda \right) 
\frac{z^{n}}{n!}=\left( \sum_{n=0}^{\infty }\mathcal{B}_{n}\left( \alpha
\mid \lambda \right) \frac{z^{n}}{n!}\right) \left( \sum_{n=0}^{\infty }x^{n}%
\frac{z^{n}}{n!}\right) \text{.}
\end{equation*}

From the rule of Cauchy product, we get%
\begin{equation}
\sum_{n=0}^{\infty }\mathcal{B}_{n}\left( x:\alpha \mid \lambda \right) 
\frac{z^{n}}{n!}=\sum_{n=0}^{\infty }\left( \sum_{k=0}^{n}\binom{n}{k}%
\mathcal{B}_{k}\left( \alpha \mid \lambda \right) x^{n-k}\right) \frac{z^{n}%
}{n!}\text{.}  \label{equ.13}
\end{equation}

Comparing the coefficients of $\frac{z^{n}}{n!}$ in (\ref{equ.13}), we have 
\begin{equation}
\mathcal{B}_{n}\left( x:\alpha \mid \lambda \right) =\sum_{k=0}^{n}\binom{n}{%
k}\mathcal{B}_{k}\left( \alpha \mid \lambda \right) x^{n-k}\text{.}
\label{equ.15}
\end{equation}

Similarly, we can get identities for Euler-type and Genocchi-type
polynomials. Therefore, we discover the following theorem.

\begin{theorem}
The following identities hold true:%
\begin{eqnarray*}
\mathcal{B}_{n}\left( x:\alpha \mid \lambda \right) &=&\sum_{k=0}^{n}\binom{n%
}{k}\mathcal{B}_{k}\left( \alpha \mid \lambda \right) x^{n-k}\text{,} \\
\mathcal{E}_{n}\left( x:\alpha \mid \lambda \right) &=&\sum_{k=0}^{n}\binom{n%
}{k}\mathcal{E}_{k}\left( \alpha \mid \lambda \right) x^{n-k}\text{,} \\
\mathcal{G}_{n}\left( x:\alpha \mid \lambda \right) &=&\sum_{k=0}^{n}\binom{n%
}{k}\mathcal{G}_{k}\left( \alpha \mid \lambda \right) x^{n-k}\text{.}
\end{eqnarray*}
\end{theorem}

Let us now apply the familiar derivative $\frac{d}{dx}$ in the both sides of
Definition \ref{Definition 1},%
\begin{eqnarray*}
\frac{d}{dx}\left( \sum_{n=0}^{\infty }\mathcal{B}_{n}\left( x:\alpha \mid
\lambda \right) \frac{z^{n}}{n!}\right) &=&\sum_{n=0}^{\infty }\left( \frac{d%
}{dx}\mathcal{B}_{n}\left( x:\alpha \mid \lambda \right) \right) \frac{z^{n}%
}{n!} \\
&=&\frac{z}{\lambda E_{\alpha }\left( z\right) -1}\left[ \frac{d}{dx}e^{xz}%
\right] \\
&=&\frac{z^{2}}{\lambda E_{\alpha }\left( z\right) -1}e^{xz} \\
&=&\sum_{n=0}^{\infty }\mathcal{B}_{n}\left( \alpha \mid \lambda \right) 
\frac{z^{n+1}}{n!}\text{.}
\end{eqnarray*}

Similarly, we can procure the derivatives of Euler-type and Genocchi-type
polynomials. Thus, we state the following theorem.

\begin{theorem}
\label{Theorem 2}The following identities hold true:%
\begin{equation*}
\frac{d}{dx}\mathcal{B}_{n}\left( x:\alpha \mid \lambda \right) =n\mathcal{B}%
_{n-1}\left( x:\alpha \mid \lambda \right) \text{, }\frac{d}{dx}\mathcal{E}%
_{n}\left( x:\alpha \mid \lambda \right) =n\mathcal{E}_{n-1}\left( x:\alpha
\mid \lambda \right)
\end{equation*}%
and 
\begin{equation*}
\frac{d}{dx}\mathcal{G}_{n}\left( x:\alpha \mid \lambda \right) =n\mathcal{G}%
_{n-1}\left( x:\alpha \mid \lambda \right) \text{.}
\end{equation*}
\end{theorem}

Polynomials $\mathcal{A}_{n}\left( x\right) $ are called Appell polynomials
when they have the following identity:%
\begin{equation}
\frac{d}{dx}\mathcal{A}_{n}\left( x:\alpha \right) =n\mathcal{A}_{n-1}\left(
x:\alpha \right) \text{.}  \label{equ.14}
\end{equation}

So, by Theorem \ref{Theorem 2} and (\ref{equ.14}), we have the following
corollary.

\begin{corollary}
Our polynomials, which are $\mathcal{B}_{n}\left( x:\alpha \mid \lambda
\right) $, $\mathcal{E}_{n}\left( x:\alpha \mid \lambda \right) $ and $%
\mathcal{G}_{n}\left( x:\alpha \mid \lambda \right) $, are Appell
polynomials.
\end{corollary}

From Theorem \ref{Theorem 2}, we have%
\begin{equation*}
\int_{0}^{1}\mathcal{B}_{n}\left( x:\alpha \mid \lambda \right) dx=\frac{%
\mathcal{B}_{n+1}\left( 1:\alpha \mid \lambda \right) -\mathcal{B}_{n}\left(
x:\alpha \mid \lambda \right) }{n+1}\text{.}
\end{equation*}

More generally,%
\begin{equation*}
\int_{x}^{x+1}\mathcal{B}_{n}\left( y:\alpha \mid \lambda \right) dy=\frac{%
\mathcal{B}_{n+1}\left( x+1:\alpha \mid \lambda \right) -\mathcal{B}%
_{n}\left( x:\alpha \mid \lambda \right) }{n+1}\text{.}
\end{equation*}

Thus, we get the following theorem.

\begin{theorem}
The following identities hold true:%
\begin{eqnarray*}
\int_{x}^{x+1}\mathcal{B}_{n}\left( y:\alpha \mid \lambda \right) dy &=&%
\frac{\mathcal{B}_{n+1}\left( x+1:\alpha \mid \lambda \right) -\mathcal{B}%
_{n}\left( x:\alpha \mid \lambda \right) }{n+1}\text{,} \\
\int_{x}^{x+1}\mathcal{E}_{n}\left( y:\alpha \mid \lambda \right) dy &=&%
\frac{\mathcal{E}_{n+1}\left( x+1:\alpha \mid \lambda \right) -\mathcal{E}%
_{n}\left( x:\alpha \mid \lambda \right) }{n+1}
\end{eqnarray*}%
and%
\begin{equation*}
\int_{x}^{x+1}\mathcal{G}_{n}\left( y:\alpha \mid \lambda \right) dy=\frac{%
\mathcal{G}_{n+1}\left( x+1:\alpha \mid \lambda \right) -\mathcal{G}%
_{n}\left( x:\alpha \mid \lambda \right) }{n+1}\text{.}
\end{equation*}
\end{theorem}

\section{\textbf{Identities including special polynomials arising from point
of view of Fractional calculus}}

Recently works involving the integral of the product of several type
Bernstein polynomials \cite{Acikgoz}, $p$-adic integral representation for $%
q $-Bernoulli numbers and studied on the $q$-integral representation of $q$%
-Bernstein polynomials \cite{Kim 4}, fermionic $p$-adic integral
representation for Frobenius-Euler numbers and polynomials \cite{Araci 6},
derivative representations of Bernstein polynomials \cite{Cetin} have been
investigated.

In this final part, we derive some new interesting related to special
polynomials by utilizing from Caputo-fractional derivative and
Riemann-Liouville integral.

As well known, Apostol-Bernoulli polynomials are given to be:%
\begin{equation}
\mathcal{F}\left( t,z\mid \lambda \right) =\frac{z}{\lambda e^{z}-1}%
e^{tz}=\sum_{n=0}^{\infty }B_{n}\left( t\mid \lambda \right) \frac{z^{n}}{n!}%
\text{.}  \label{equ.16}
\end{equation}

Note that Apostol-Bernoulli polynomials are analytic on the region $%
\mathfrak{D}=\left\{ z\in 
\mathbb{C}
\mid \left\vert z+\log \lambda \right\vert <2\pi \right\} $ (see \cite{Luo1}%
, \cite{He1}, \cite{He2}).

Differentiating in the both sides of (\ref{equ.16}), we have%
\begin{equation}
\frac{d}{dt}B_{n}\left( t\mid \lambda \right) =nB_{n-1}\left( t\mid \lambda
\right) \text{ (see \cite{Luo1}, \cite{He1}, \cite{He2}).}  \label{equ.17}
\end{equation}

When $t=0$ in (\ref{equ.16}), we have $B_{n}\left( 0\mid \lambda \right)
:=B_{n}\left( \lambda \right) $ are called Bernoulli numbers, which can be
generated by%
\begin{equation}
\mathcal{F}\left( z\mid \lambda \right) =\frac{z}{\lambda e^{z}-1}%
=\sum_{n=0}^{\infty }B_{n}\left( \lambda \right) \frac{z^{n}}{n!}.
\label{equ.18}
\end{equation}

By (\ref{equ.16}) and (\ref{equ.18}), we have the following functional
equation:%
\begin{equation*}
\mathcal{F}\left( t,z\mid \lambda \right) =e^{tz}\mathcal{F}\left( z\mid
\lambda \right) .
\end{equation*}

By using Taylor's formula in the last identity, we have%
\begin{equation*}
B_{m}\left( t\mid \lambda \right) =\sum_{k=0}^{m}\binom{m}{k}%
t^{m-k}B_{k}\left( \lambda \right) =\sum_{k=0}^{m}\binom{m}{k}%
t^{k}B_{m-k}\left( \lambda \right) \text{ (see \cite{Luo1}, \cite{He1}, \cite%
{He2}).}
\end{equation*}

Let us now take $f\left( t\right) =B_{m}\left( t\mid \lambda \right) $ in (%
\ref{equ.7}), leads to%
\begin{eqnarray*}
\mathcal{D}^{\alpha }B_{m}\left( t\mid \lambda \right) &=&\frac{1}{\Gamma
\left( n-\alpha \right) }\int_{0}^{t}\frac{\frac{d^{n}}{dt^{n}}B_{m}\left(
t\mid \lambda \right) \mid _{t=s}}{\left( t-s\right) ^{\alpha -n+1}}ds \\
&=&m\left( m-1\right) \cdots \left( m-n+1\right) \sum_{k=0}^{m-n}\binom{m-n}{%
k}B_{m-n-k}\left( \lambda \right) \\
&&\times \left[ \frac{1}{\Gamma \left( n-\alpha \right) }\int_{0}^{t}\frac{%
s^{k}}{\left( t-s\right) ^{\alpha -n+1}}ds\right] \\
&=&\frac{\Gamma \left( m+1\right) }{\Gamma \left( m-n+1\right) }%
\sum_{k=0}^{m-n}\frac{k!\binom{m-n}{k}B_{m-n-k}\left( \lambda \right) }{%
\Gamma \left( n+k-\alpha +1\right) }t^{k-\alpha +n}\text{.}
\end{eqnarray*}

Therefore, we procure the following theorem.

\begin{theorem}
The following identity%
\begin{equation*}
\mathcal{D}^{\alpha }B_{m}\left( t\mid \lambda \right) =\frac{\Gamma \left(
m+1\right) }{\Gamma \left( m-n+1\right) }\sum_{k=0}^{m-n}\frac{k!\binom{m-n}{%
k}B_{m-n-k}\left( \lambda \right) }{\Gamma \left( n+k-\alpha +1\right) }%
t^{k-\alpha +n}
\end{equation*}%
is true.
\end{theorem}

In \cite{Luo1}, Apostol-Bernoulli polynomials of higher order are defined by%
\begin{equation}
\underset{h-times}{\underbrace{\frac{z}{\lambda e^{z}-1}\frac{z}{\lambda
e^{z}-1}\cdots \frac{z}{\lambda e^{z}-1}}}e^{tz}=\sum_{m=0}^{\infty
}B_{m}^{\left( h\right) }\left( t\mid \lambda \right) \frac{z^{n}}{n!}\text{.%
}  \label{equ.19}
\end{equation}

Note that $B_{m}^{\left( h\right) }\left( t\right) $ is analytic on $%
\mathfrak{D}$. It follows from (\ref{equ.19}), we have 
\begin{equation}
\frac{d}{dt}B_{m}^{\left( h\right) }\left( t\mid \lambda \right)
=mB_{m-1}^{\left( h\right) }\left( t\mid \lambda \right) \text{ \textit{and} 
}\frac{d^{n}}{dt^{n}}B_{m}^{\left( h\right) }\left( t\mid \lambda \right) =%
\frac{\Gamma \left( m+1\right) }{\Gamma \left( m-n+1\right) }B_{m-n}^{\left(
h\right) }\left( t\mid \lambda \right) \text{ (see \cite{Luo1}).}
\label{equ.20}
\end{equation}

Substituting $x=0$ into (\ref{equ.19}), $B_{m}^{\left( h\right) }\left(
0\mid \lambda \right) :=B_{m}^{\left( h\right) }\left( \lambda \right) $ are
called Bernoulli polynomials of higher order.

Owing to (\ref{equ.19}) and (\ref{equ.20}), we readily see that%
\begin{eqnarray*}
\mathcal{D}^{\alpha }B_{m}^{\left( h\right) }\left( t\mid \lambda \right) &=&%
\frac{1}{\Gamma \left( n-\alpha \right) }\int_{0}^{t}\frac{\frac{d^{n}}{%
dt^{n}}B_{m}^{\left( h\right) }\left( t\mid \lambda \right) \mid _{t=s}}{%
\left( t-s\right) ^{\alpha -n+1}}ds \\
&=&m\left( m-1\right) \cdots \left( m-n+1\right) \sum_{k=0}^{m-n}\binom{m-n}{%
k}B_{m-n-k}^{\left( h\right) }\left( \lambda \right) \\
&&\times \left[ \frac{1}{\Gamma \left( n-\alpha \right) }\int_{0}^{t}\frac{%
s^{k}}{\left( t-s\right) ^{\alpha -n+1}}ds\right] \\
&=&\frac{\Gamma \left( m+1\right) }{\Gamma \left( m-n+1\right) }%
\sum_{k=0}^{m-n}\frac{k!\binom{m-n}{k}B_{m-n-k}^{\left( h\right) }\left(
\lambda \right) }{\Gamma \left( n+k-\alpha +1\right) }t^{k-\alpha +n} \\
&=&\frac{\Gamma \left( m+1\right) }{\Gamma \left( m-n+1\right) }%
\sum_{k=0}^{m-n}\frac{k!\binom{m-n}{k}}{\Gamma \left( n+k-\alpha +1\right) }%
t^{k-\alpha +n} \\
&&\times \left( \sum_{\underset{s_{h}\geq 0}{s_{1}+s_{2}+\cdots +s_{h}=m-n-k}%
}\binom{m-n-k}{s_{1},s_{2},\cdots ,s_{h}}\left(
\tprod\limits_{j=1}^{h}B_{s_{j}}\right) \right) \text{.}
\end{eqnarray*}

Therefore, we can state the following theorem.

\begin{theorem}
The following equality%
\begin{eqnarray*}
\mathcal{D}^{\alpha }B_{m}^{\left( h\right) }\left( t\mid \lambda \right) &=&%
\frac{\Gamma \left( m+1\right) }{\Gamma \left( m-n+1\right) }\sum_{k=0}^{m-n}%
\frac{k!\binom{m-n}{k}}{\Gamma \left( n+k-\alpha +1\right) }t^{k-\alpha +n}
\\
&&\times \left( \sum_{\underset{s_{h}\geq 0}{s_{1}+s_{2}+\cdots +s_{h}=m-n-k}%
}\binom{m-n-k}{s_{1},s_{2},\cdots ,s_{h}}\left(
\tprod\limits_{j=1}^{h}B_{s_{j}}\right) \right)
\end{eqnarray*}%
is true.
\end{theorem}

We recall the definition of generating function for Bernoulli-type
polynomials as follows:%
\begin{equation*}
\sum_{n=0}^{\infty }\mathcal{B}_{n}\left( x:\alpha \mid \lambda \right) 
\frac{z^{n}}{n!}=\frac{z}{\lambda E_{\alpha }\left( z\right) -1}e^{xz}
\end{equation*}%
and also%
\begin{equation*}
\frac{d}{dx}\mathcal{B}_{n}\left( x:\alpha \mid \lambda \right) =n\mathcal{B}%
_{n-1}\left( x:\alpha \mid \lambda \right) \text{.}
\end{equation*}

Taking $f\left( t\right) =\mathcal{B}_{n}\left( x:\alpha \mid \lambda
\right) $ in (\ref{equ.7}), we compute%
\begin{eqnarray*}
D^{\alpha }\mathcal{B}_{n}\left( x:\alpha \mid \lambda \right)  &=&\frac{1}{%
\Gamma \left( n-\alpha \right) }\int_{0}^{t}\frac{\frac{d^{n}}{dt^{n}}%
\mathcal{B}_{n}\left( x:\alpha \mid \lambda \right) \mid _{t=s}}{\left(
t-s\right) ^{\alpha -n+1}}ds \\
&=&m\left( m-1\right) \cdots \left( m-n+1\right) \sum_{k=0}^{m-n}\binom{m-n}{%
k}\mathcal{B}_{n}\left( \alpha \mid \lambda \right)  \\
&&\times \left[ \frac{1}{\Gamma \left( n-\alpha \right) }\int_{0}^{t}\frac{%
s^{k}}{\left( t-s\right) ^{\alpha -n+1}}ds\right]  \\
&=&\frac{\Gamma \left( m+1\right) }{\Gamma \left( m-n+1\right) }%
\sum_{k=0}^{m-n}\frac{k!\binom{m-n}{k}\mathcal{B}_{n}\left( \alpha \mid
\lambda \right) }{\Gamma \left( n+k-\alpha +1\right) }t^{k-\alpha +n}\text{.}
\end{eqnarray*}

We can acquire similar identities for Euler-type polynomials and
Genocchi-type polynomials. So, we state the following interesting theorem.

\begin{theorem}
The following equalities hold true:%
\begin{eqnarray*}
D^{\alpha }\mathcal{B}_{n}\left( x:\alpha \mid \lambda \right) &=&\frac{%
\Gamma \left( m+1\right) }{\Gamma \left( m-n+1\right) }\sum_{k=0}^{m-n}\frac{%
k!\binom{m-n}{k}\mathcal{B}_{n}\left( \alpha \mid \lambda \right) }{\Gamma
\left( n+k-\alpha +1\right) }t^{k-\alpha +n}\text{,} \\
D^{\alpha }\mathcal{E}_{n}\left( x:\alpha \mid \lambda \right) &=&\frac{%
\Gamma \left( m+1\right) }{\Gamma \left( m-n+1\right) }\sum_{k=0}^{m-n}\frac{%
k!\binom{m-n}{k}\mathcal{E}_{n}\left( \alpha \mid \lambda \right) }{\Gamma
\left( n+k-\alpha +1\right) }t^{k-\alpha +n}\text{,}
\end{eqnarray*}%
and%
\begin{equation*}
D^{\alpha }\mathcal{G}_{n}\left( x:\alpha \mid \lambda \right) =\frac{\Gamma
\left( m+1\right) }{\Gamma \left( m-n+1\right) }\sum_{k=0}^{m-n}\frac{k!%
\binom{m-n}{k}\mathcal{G}_{n}\left( \alpha \mid \lambda \right) }{\Gamma
\left( n+k-\alpha +1\right) }t^{k-\alpha +n}\text{.}
\end{equation*}
\end{theorem}

\end{document}